\documentclass[reqno]{amsart}
\usepackage{amsmath}
\usepackage{amscd}
\usepackage{graphics}
\usepackage{latexsym}
\usepackage{hyperref}

\theoremstyle{plain}
\newtheorem{theorem}{Theorem}[section]
\newtheorem{thm}[theorem]{Theorem}
\newtheorem{cor}[theorem]{Corollary}
\newtheorem{lem}[theorem]{Lemma}
\newtheorem{prop}[theorem]{Proposition}

\theoremstyle{definition}
\newtheorem{defn}[theorem]{Definition}

\theoremstyle{remark}

\newcommand{\QQ}{\mathbb{Q}}

\newcommand{\AAA}{\mathbb{A}}
\newcommand{\PP}{\mathbb{P}}

\newcommand{\SP}{\text{Spec }}

\newcommand{\marpar}[1]{}

\newcommand{\lt}{\left}
\newcommand{\rt}{\right}

\newcommand{\OO}{\mathcal O}

\def\PP{{\mathbb P}}

\begin{document}

\title[Free rational quotient]
{The maximal free rational quotient}

\author[Starr]{Jason Michael Starr}
\address{Department of Mathematics \\
  Massachusetts Institute of Technology \\ Cambridge MA 02139}
\email{jstarr@math.mit.edu} 

\date{\today}

\begin{abstract}
This short, expository note proves existence of the maximal quotient
of a variety by free rational curves.  
\end{abstract}

\maketitle


\section{Definition of a maximal free rational quotient}~\label{sec-sor}
\marpar{sec-sor}

\begin{defn} \label{defn-0}
Let $V$ be a Deligne-Mumford stack over a field $k$, and denote the
smooth locus by $V^\text{sm} \subset V$.
A $1$-morphism $f:\PP^1_k \rightarrow V^\text{sm}$ is a 
\emph{free rational curve to $V$} 
if $f^*T_V$ is generated by global
sections and has positive degree.
\end{defn}

\medskip\noindent
Let $S$ be an irreducible algebraic space, 
let $\pi_{\overline{X}}:\overline{X}\rightarrow S$ be a proper,
locally finitely presented $1$-morphism of Deligne-Mumford stacks
with integral geometric generic fiber, and let $X\subset
\overline{X}$ be a normal dense open substack.  Denote by $\pi_X:X
\rightarrow S$ the restriction of $\pi_{\overline{X}}$.

\begin{defn} \label{defn-1}
\marpar{defn-1}    
A \emph{free rational
  quotient} of $\pi_X$ is a triple $(X^*,Q^*,\phi)$
where $X^*\subset X$ is a dense open substack, where $Q^*$ is a normal
algebraic space, finitely presented over $S$ with integral
geometric generic fiber,
and where  
$\phi:X^* \rightarrow Q^*$ 
is a dominant $1$-morphism of $S$-stacks satisfying,
\begin{enumerate}
\item[(i)] the geometric generic fiber $F$ of $\phi$ is integral, and
\item[(ii)] a general pair of distinct points of $F$ is contained in
  the image 
  of a free rational curve.
\end{enumerate}
A free rational quotient is \emph{trivial} if $\phi:X^* \rightarrow
Q^*$ is birational, and \emph{nontrivial} otherwise.  

\noindent
A free rational quotient $(X^*,Q^*,\phi:X^*\rightarrow Q^*)$ is 
\emph{maximal} if 
for every free
rational quotient $(X^*_1,Q^*_1,\phi_1:X^*_1\rightarrow Q^*_1)$ there
exists a dense open subset $U\subset Q^*_1$ and a smooth morphism
$\psi:U \rightarrow Q^*$ such that $\phi|_{\phi_1^{-1}(U)} = \psi\circ
\phi_1$.  
\end{defn}

\begin{thm} \label{thm-2}
\marpar{thm-2}
There exists a maximal free rational quotient.
\end{thm}

\medskip\noindent
It is not true that a maximal free rational quotient is unique, but it
is unique up to unique birational equivalence of $Q^*$.  

\section{Proof of Theorem \ref{thm-2}}\label{sec-proofs}
\marpar{sec-proofs}
\noindent
The proof is very similar to the proofs of existence of the
rational quotient in ~\cite{De} and ~\cite{K}.  Existence of
the maximal free rational quotient can be deduced from theorems there.
However there are $2$ special features of this case: The relation of
containment in a free rational curve is already a rational equivalence
relation, so existence of the quotient is less technical than the
general case.  
Also, unlike the general case, there is no need to perform a purely
inseparable base-change of $S$ to define the quotient.  Of
course if $\OO_S$ contains $\QQ$ and if $\pi$ is smooth and proper,
~\cite[1.1]{KMM92c} implies the free rational quotient is the rational
quotient.  

\medskip\noindent
If there is no free rational curve in the geometric generic fiber of
$\pi_X$, then the trivial
rational quotient $(X,X,\text{Id}_X:X\rightarrow X)$ is a maximal free
rational quotient.  Therefore assume there is a free rational
curve to $X$ whose image is contained in the smooth locus of $\pi_X$
(equivalently, there is a free rational curve in the geometric
generic fiber of $\pi_X$).
Use will be made of the flat, proper morphism $C\rightarrow S$
obtained from $2$ copies of $\PP^1_S$ by identifying the $0$ section
in the first copy to the $0$ section in the second copy.

\medskip\noindent
Denote by
$\text{Hom}_S(\PP^1_S, X)$ the Deligne-Mumford stack constructed in
~\cite{OHom}.  Define $H \subset \text{Hom}_S(\PP^1_S,X)$ to be the open
substack parametrizing free rational curves in fibers of the smooth
locus of $\pi_X$.
By hypothesis, $H$ is nonempty. 

\medskip\noindent
Let
$H_i \subset H$ be a connected component.  Denote by,
$$
u_i: H_i \times_S \PP^1_S \rightarrow \overline{X},
$$
and by,
$$
u_i^{(2)}: H_i \times_S \PP^1_S \times_S \PP^1_S \rightarrow \overline{X}
\times_S \overline{X},
$$
the obvious $1$-morphisms.

\begin{lem} \label{lem-3}
\marpar{lem-3}
The $1$-morphism $u_i$ is smooth.
\end{lem} 

\begin{proof}
The proof is the same as in the case that $\pi_X$ is a projective
morphism, ~\cite[Cor. II.3.5.4]{K}. 
\end{proof}

\medskip\noindent
Denote by $W_i \subset \overline{X} \times_S \overline{X}$ the closed
image of 
$u_i^{(2)}$, i.e., $W_i$ is the minimal closed substack such that
$u_i^{(2)}$ factors through $W_i$.  For any geometric point $x$ of
$\overline{X}$ with residue field $\kappa(x)$, 
denote by $(W_i)_x \subset X\otimes_{\OO_S} \kappa(x)$ 
the scheme
$\text{pr}_2(\text{pr}_1^{-1}(x) \cap W_i)$.  

\medskip\noindent
Because $H_i\times_S \PP^1_S
\times_S \PP^1_S$ is irreducible and reduced, also $W_i$ is
irreducible and reduced.  Because $u_i$ is smooth, $\text{pr}_1:W_i
\rightarrow \overline{X}$ is surjective on geometric points.  The geometric
generic fiber has pure dimension 
$$
d_i = \text{dim}(W_i\otimes_{\OO_S}
K(S)) - \text{dim}(X\otimes_{\OO_S} K(S)).
$$
Of course $d_i$ is bounded by $\text{dim}(X\otimes_{\OO_S} K(S))$.  
Let $H_i$ be a connected component such that $d_i$ is maximal.

\medskip\noindent
Let $H_j$ be any connected component of $H$.  
Denote by $V_{i,j}$ the
2-fiber product,
$$
V_{i,j} = (H_i \times_S \PP^1_S) \times_{u_i,X,u_j} (H_j \times_S
\PP^1_S).
$$
In other words $V_{i,j}$ parametrizes data
$(([f_i],t_i),([f_j],t_j),\theta)$ where $f_i$, resp. $f_j$, is an object of
$H_i$, resp. $H_j$, where $t_i,t_j$ are points of $\PP^1$, and
where
$\theta:f_i(t_i) \rightarrow f_j(t_j)$ is an equivalence of objects. 
Denote by $F_{i,j}$ the $1$-morphism of $S$-stacks,
$$
F_{i,j}: V_{i,j}\times_S \PP^1_S \times_S \PP^1_S \rightarrow
\overline{X}\times_S \overline{X},
$$
that sends a datum $((([f_i],t_i),([f_j],t_j),\theta), t'_i,t'_j)$ to
$(f_i(t_i'),f_j(t_j'))$.  
Alternatively $F_{i,j}$ is the 1-morphism whose domain is,
$$
(H_i \times_S \PP^1_S \times_S \PP^1_2)
\times_{u_i\circ
  \text{pr}_{1,3},X,u_j\circ \text{pr}_{1,2}} (H_j \times_S \PP^1_S
\times_S \PP^1_S),
$$
such that $\text{pr}_1\circ F_{i,j}$ is $u_i\circ \text{pr}_{1,2}$ and
such that $\text{pr}_2\circ F_{i,j}$ is $u_j\circ \text{pr}_{1,3}$.

\begin{prop} \label{prop-4}
\marpar{prop-4}
The image of $F_{i,j}$ is contained in $W_i$.
\end{prop}

\begin{proof}
\medskip\noindent
Note that $V_{i,j}$ is smooth over $S$.  Moreover, because $u_i$ and
$u_j$ are smooth, $V_{i,j}$ is nonempty.  Let $((f_i,t_i),(f_j,t_j),\theta)$
be a point of $V_{i,j}$.  There is a reducible, connected genus $0$
curve $C$ 
obtained by identifying $t_i$ in one copy of $\PP^1$ to $t_j$ in a
second copy of $\PP^1$.  The morphisms $f_i$, $f_j$ and the
equivalence $\theta$ induce a
1-morphism $f:C\rightarrow X$ whose restriction to the
first irreducible component is $f_i$ and whose restriction to the
second irreducible component is $f_j$.  In a suitable sense,
$f:C\rightarrow X$ is still a free rational curve, and it deforms to
free rational curves $f':\PP^1 \rightarrow X$.  After examining these
deformations, the proposition easily follows.

\medskip\noindent
Let $0:S \rightarrow \AAA^1_S \times_S \PP^1_S$ be the morphism whose
projection to each factor is the zero section.  Denote by $Z \subset
\AAA^1_S \times_S \PP^1_S$ the image of $0$.  Denote by $P$ the
blowing up of $\AAA^1_S \times_S \PP^1_S$ along $Z$.  The projection
morphism $\text{pr}_{\AAA^1}:P \rightarrow \AAA^1_S$ is flat and
projective.  Moreover, the restriction of $P$ over $\mathbb{G}_{m,S}
\subset \AAA^1_S$ is canonically isomorphic to
$\mathbb{G}_{m,S}\times_S \PP^1_S$.  And the restriction of $P$ over
the zero section of $\AAA^1_S$ is canonically isomorphic to the curve
$C$
over $S$
obtained by identifying $0$ in one copy of $\PP^1$ to $0$ in a second
copy of $\PP^1$.

\medskip\noindent
Denote $S'=\AAA^1_S$.  Denote by
$\text{Hom}_{S'}(P,X_{S'})$ the Deligne-Mumford stack constructed in
~\cite{OHom}.  Denote by $H' \subset \text{Hom}_{S'}(P,X_{S'})$ the open
substack parametrizing morphisms to the smooth locus of $(\pi_X)'$
such that
the pullback of $T_{(\pi_X)'}$ restricts to
a globally generated sheaf of positive degree
on every irreducible component of every
geometric fiber.  The morphism $H' \rightarrow S'$ is smooth for 
reasons similar to ~\cite[Cor. II.3.5.4]{K}.  

\medskip\noindent
Let $V_{i,j,k}$ be a connected component of $V_{i,j}$.
Define $f_{V,i}:V_{i,j,k}\times_S \PP^1_S \rightarrow X$ to be
the composition 
$$
u_i\circ (\text{pr}_{H_i}\circ
\text{pr}_{H_i\times_S \PP^1_S},\text{Id}_{\PP^1_S}):
V_{i,j,k}\times_S \PP^1_S \rightarrow H_i \times_S \PP^1_S \rightarrow
X.
$$
Define $f_{V,j}:V_{i,j,k}\times_S \PP^1_S \rightarrow X$ to be the
composition $u_j\circ (\text{pr}_{H_j}\circ \text{pr}_{H_j\times_S
  \PP^1_S},\text{Id}_{\PP^1_S})$.
Define $s_i:V_{i,j,k} \rightarrow V_{i,j,k} \times_S \PP^1_S$ to be
the unique $V_{i,j,k}$-morphism such that $\text{pr}_{\PP^1_S}\circ
s_i = \text{pr}_{\PP^1_S} \circ \text{pr}_{H_i \times_S \PP^1_S}$.
Define $s_j:V_{i,j,k} \rightarrow V_{i,j,k} \times_S \PP^1_S$
similarly.

\medskip\noindent
Replacing $V_{i,j,k}$ by a dense open subset,
there exist $2$ isomorphisms of $V_{i,j,k}$-schemes, 
$$
\alpha_i, \alpha_j :V_{i,j,k} \times_S \PP^1_S \rightarrow V_{i,j,k}
\times_S \PP^1_S,
$$
such that $s_i = \alpha_i\circ 0$ and $s_j = \alpha_j\circ 0$ 
where $0$ is the zero section of $V_{i,j,k}\times_S \PP^1_S
\rightarrow V_{i,j,k}$.  There is a unique $1$-morphism of $S$-stacks,
$$
f_{i,j,k}:V_{i,j,k} \times_S C \rightarrow X,
$$
such that the restriction of $f_{i,j,k}$ to the first irreducible
component of $V_{i,j,k}\times_S C$ is $f_{V,i}\circ \alpha_i$, and the
restriction to the second irreducible component is $f_{V,j}\circ
\alpha_j$.  The image of $f_{i,j,k}$ is contained in
the smooth locus of $\pi_X$, and the restriction of
$f_{i,j,k}^*T_{\pi_X}$ to each irreducible component is generated by
global sections relative to $V_{i,j,k}$.  Denote by,
$$
f_{i,j,k}^{(2)}: V_{i,j,k} \times_S C \times_S C \rightarrow
\overline{X}\times_S \overline{X},
$$
the obvious $1$-morphism.

\medskip\noindent
By definition of $H'$
there is a
$1$-morphism, 
$$
q:V_{i,j,k} \rightarrow S\times_{0,S'} H',
$$ 
such that the pullback by $q$ of the universal morphism is
$2$-equivalent to
$f_{i,j,k}$.  Because $V_{i,j,k}$ is connected, the image of $q$ is
contained in a connected component $H'_l$ of $H'$.  By definition of
$H$ there is a $1$-morphism of $S$-stacks,
$$
r: \mathbb{G}_{m,S} \times_{S'} H'_l \rightarrow H,
$$
such that the restriction to $\mathbb{G}_{m,S}\times_{S'} H'_l$ 
of the universal morphism over $H'_l$
is $2$-equivalent to the
pullback by $r$ of the universal morphism over $H$.  
Denote by $H_l$ the
connected component of $H$ dominated by $r$.  The morphism $r$
dominates a connected component $H_l \subset H$.
There exists a 
$1$-isomorphism, 
$$
i:\mathbb{G}_{m,S} \times_{S'} H'_l \rightarrow \mathbb{G}_{m,S}
\times_S H_l,
$$
unique up to unique $2$-equivalence, 
such that $\text{pr}_{\mathbb{G}_m}\circ i =
\text{pr}_{\mathbb{G}_m}$ and such that $\text{pr}_{H_l} \circ i$ is
$2$-equivalent to $r$.

\medskip\noindent
Denote by,
$$
v_l: H'_l \times_{S'} P \rightarrow \overline{X},
$$
and by,
$$
v_l^{(2)}:H'_l \times_{S'} P \times_{S'} P \rightarrow
\overline{X}\times_S \overline{X}, 
$$
the obvious morphisms.  Denote by $W'_l \subset \overline{X}\times_S
\overline{X}$ the
minimal closed substack through which $v_l^{(2)}$ factors.  Because
$\text{pr}_{S'}: H'_l \times_{S'} P \times_{S'} P \rightarrow S'$ is
flat, 
the preimage of $\mathbb{G}_{m,S} \subset S'$ is dense.
Thus $W'_l$ is the image of the restriction of $v_l^{(2)}$ over
$\mathbb{G}_{m,S}$.  The restriction of $v_l^{(2)}$ is $2$-equivalent
to $u_l^{(2)}\circ (r,\text{Id}_{\PP^1},\text{Id}_{\PP^1})$.
Therefore the image of $v_l^{(2)}$ equals the image of $u_l^{(2)}$,
i.e., $W'_l = W_l$.  

\medskip\noindent
On the other hand, the pullback of $v_l^{(2)}$ to
$V_{i,j,k} \times_S C \times_S C$ is $2$-equivalent to $f_{i,j,k}^{(2)}$.
There are $2$ irreducible components of $C$, and thus $4$ irreducible
components of $C\times_S C$.
Restrict $f_{i,j,k}^{(2)}$
to the irreducible component of $C\times_S C$
that is the product
of the first irreducible component of $C$ and the first irreducible
component of $C$.  This is $2$-equivalent to the pullback of
$u_i^{(2)}$, hence $W_l$ contains $W_i$.
Because $W_l$ is an
integral stack of dimension at most $d_i$ containing the
$d_i$-dimensional stack $W_i$,
$W_l$ equals $W_i$.  

\medskip\noindent
Finally, restrict $f_{i,j,k}^{(2)}$ to the irreducible component of
$C\times_S C$ that is the product of the first irreducible component
of $C$ and the second irreducible component of $C$.  This is
$2$-equivalent to the pullback of $F_{(i,j)}$, hence $W_i=W_l$ contains
the image of $F_{(i,j)}$.
\end{proof}

\begin{lem} \label{lem-4.5}
\marpar{lem-4.5}
The geometric generic fiber of $\text{pr}_1:W_i \rightarrow \overline{X}$ is
integral.
\end{lem}

\begin{proof}
Denote by $K$ the algebraic closure of the function field of $\overline{X}$.
Since $u_i:H_i \times_S \PP^1_S \rightarrow \overline{X}$ is smooth, the
geometric generic fiber $(H_i\times_S
\PP^1_S)\otimes_{\OO_{\overline{X}}} K$ is 
smooth over $K$.  There is an induced $1$-morphism,
$$
(u_i^{(2)})_K: (H_i\times_S \PP^1_S \otimes_{ \OO_{\overline{X}} } K)
\times_{\SP(K)} 
\PP^1_K \rightarrow (\overline{X} \times_S
\overline{X})\otimes_{\text{pr}_1,\overline{X}} \SP(K) \cong
\overline{X}\otimes_{\OO_S} K
$$
Formation of the closed image is compatible with flat
base change.  Therefore the closed image of $(u_i^{(2)})_K$ is
$W_i\otimes_{\text{pr}_1,\overline{X}} \SP(K)$.  In particular,
$W_i\otimes_{\text{pr}_1,\overline{X}} \SP(K)$ is reduced since the
closed image of a reduced stack is reduced.  

\medskip\noindent
The proof that the geometric generic fiber is irreducible is
essentially the same as the proof of Proposition~\ref{prop-4}.  Let
$W' \subset W_i \otimes_{\text{pr}_1,\overline{X}} \SP(K)$ be an irreducible
component.  To prove that $W' =  W_i
\otimes_{\text{pr}_1,\overline{X}} \SP(K)$, 
it suffices to prove that it contains the image of $(u_i^{(2)})_K$.
Let $x\in \overline{X}\otimes_{\OO_S} K$ be the $K$-point corresponding to the
diagonal; $x$ is contained in $\text{Image}(u_i^{(2)})_K$.  Let $y_1$
be a $K$-point of $W'\cap \text{Image}(u_i^{(2)})_K$ and let $y_2$ be a
$K$-point of $\text{Image}(u_i^{(2)})_K$.  There are free
$K$-morphisms, $f_1, f_2:\PP^1_K \rightarrow X\otimes_{\OO_S} K$ such
that $f_1(0) = f_2(0) = x$ and $f_1(\infty)=y_1, f_2(\infty)=y_2$.
This defines a morphism from $C\otimes_{\OO_S} K$ to $X\otimes_{\OO_S}
K$.  As in the proof of Proposition~\ref{prop-4}, deformations of this
morphism are free rational curves that come from a connected component
$H_l$ of $H$.  By construction, there is an irreducible component of 
$W_l\otimes_{\text{pr}_1,X} \SP(K)$ that contains $W$ and $y_2$.
Since the dimension of $W_l$ is at most $d_i$, this irreducible
component equals $W$.  Therefore $y_2 \in W$, i.e., $W = W_i
\otimes_{\text{pr}_1,X} \SP(K)$.
\end{proof}

\medskip\noindent
Consider the projection $\text{pr}_1:W_i \rightarrow \overline{X}$.  By
~\cite[Thm. 3.2]{OS}, there exists a dense open subset 
$X^\text{flat}\subset X$
over which $W_i$ is flat.  
Denote $W_i^\text{flat} = W_i
\times_{\text{pr}_1,X} X^\text{flat}$.  
By Lemma~\ref{lem-4.5}, there is a dense
open substack $X^0 \subset X^\text{flat}$ such that every geometric
fiber of $W_i \times_{\text{pr}_1,X} X^0 \rightarrow X^0$ is
integral.  Denote $W_i^0 = W_i \times_{\text{pr}_1,\overline{X}} X^0$.

\medskip\noindent
Let $H_j \subset H$ be a connected
component and denote by $G_{i,j}$ the unique $1$-morphism,
$$
G_{i,j}: (H_j \times_S \PP^1_S \times_S \PP^1_S) \times_{u_j\circ
  \text{pr}_{1,3}, \overline{X} ,\text{pr}_1} W_i \rightarrow
\overline{X}\times_S \overline{X},
$$
such that, 
$$
\text{pr}_1\circ G_{i,j} = u_j \circ \text{pr}_{1,2} \circ
\text{pr}_{H_j\times \PP^1 \times \PP^1},
$$
and such that,
$$
\text{pr}_2\circ G_{i,j} = \text{pr}_2 \circ \text{pr}_{W_i}.
$$

\begin{cor}\label{cor-5}
\marpar{cor-5}
\begin{enumerate}
\item[(i)] The image of $G_{i,j}$ is contained in $W_i$.
\item[(ii)] For every geometric point $s\in S$, for every free
  morphism $f:\PP^1_s \rightarrow X_s$ and for every point $x\in X_s$
  such that $f(\PP^1) \cap (W_i)_x$ is nonempty, $f(\PP^1)$ is
  contained in $(W_i)_x$.
\item[(iii)] The image of the ``composition morphism'',
$$
c:W_i \times_{\text{pr}_2,\overline{X},\text{pr}_1} W_i^0 \rightarrow
\overline{X}\times_S \overline{X},
$$
is contained in $W_i$.
\item[(iv)] For every geometric point $s\in S$, for every pair of
  closed points $(x,y)
  \in \overline{X}_s \times X^0_s$, if 
  $y\in (W_i)_x$ then $(W_i)_y \subset
  (W_i)_x$.
\item[(v)] For every geometric point $s\in S$, for every pair of
  closed points $(x,y) \in X_s^0 \times X_s^0$, $y_i \in (W_i)_x$ iff
  $x\in (W_i)_y$ iff $(W_i)_x = (W_i)_y$.  
\end{enumerate}
\end{cor}

\begin{proof}
\textbf{(i):}  First of all, the projection morphism,
$$
\text{pr}_{W_i}: (H_j\times_S \PP^1_S \times_S \PP^1_S)
\times_{\overline{X}} W_i \rightarrow W_i,
$$
is smooth.  Therefore every connected component of the domain is
integral and 
dominates $W_i$.  So to prove $W_i$ contains the image of $G_{i,j}$ 
it suffices to first base-change by,
$$
u_i^{(2)}: H_i \times_S \PP^1_S \times_S \PP^1_S \rightarrow W_i.
$$
The base-change of $G_{i,j}$ by $u_i^{(2)}$ is $2$-equivalent to
$F_{i,j}$.  By Proposition~\ref{prop-4} the $W_i$ contains the
image of $F_{i,j}$. 
Therefore $W_i$ contains the image of $G_{i,j}$.

\medskip\noindent
\textbf{(ii):}  By construction, $W_i \subset \overline{X}\times_S
\overline{X}$ is 
symmetric with respect to permuting the factors.
Let $t'\in \PP^1$ be a point such that $x'=f(t')$ is in
$(W_i)_x$.  Let $H_j$ be the connected component of $H$ that contains $[f]$.
Then the subset,
$$
\lt\{ (([f],t,t'),(x',x)) \in (H_j\times_S \PP^1_S\times_S \PP^1_S)
\times W_i | t\in \PP^1_s \rt\},
$$ 
is contained in,
$$
(H_j \times_S \PP^1_S \times_S \PP^1_S)\times_{\overline{X}} W_i.
$$
Therefore by (i), $W_i$ contains the image under $G_{i,j}$.
Because $W_i$ is symmetric, this implies that $(W_i)_x$ contains
$f(\PP^1)$,

\medskip\noindent
\textbf{(iii):} The $1$-morphism $c$ satisfies
$\text{pr}_1\circ c = \text{pr}_1\circ \text{pr}_1$ and 
$\text{pr}_2\circ c = \text{pr}_2 \circ \text{pr}_2$, up to
$2$-equivalence.  
Because
$\text{pr}_1:W_i^\text{flat} \rightarrow \overline{X}$ is flat and the
geometric 
fibers are integral also the projection,
$$
\text{pr}_1: W_i \times_{\text{pr}_2,\overline{X},\text{pr}_1} W_i^\text{flat}
\rightarrow W_i,
$$
is flat and the geometric fibers are integral.  In particular the
domain is integral.  
Hence to prove $W_i$ contains the image of $c$
it suffices to first base-change by,
$$
u_i^{(2)}: H_i\times_S \PP^1_S \times_S \PP^1_S \rightarrow W_i.
$$
After base-change, this morphism factors through $G_{i,i}$.  By
(i), the $W_i$ contains the image of $G_{i,i}$.  Therefore $W_i$
contains the image of $c$.

\medskip\noindent
\textbf{(iv) and (v):} Item (iv) follows immediately from (iii), and
Item (v) follows from (iv) and symmetry of $W_i$.
\end{proof}

\begin{lem} \label{lem-7}
\marpar{lem-7}
Let $k$ be a field, let $g:Y\rightarrow Z$ be a morphism of smooth
Deligne-Mumford stacks over $k$, and let 
$f:\PP^1_k \rightarrow Y$ be a free morphism
such that $f(\PP^1)$ is contained in a fiber of $g$.
Denote by $Y^\text{sm}$ the smooth locus of $g$.  If $f(\PP^1) \cap
Y^\text{sm}$ is nonempty, then $f(\PP^1) \subset Y^\text{sm}$.
\end{lem}

\begin{proof}
There is a morphism of locally free sheaves on $\PP^1_k$,
$$
dg:f^* T_Y \rightarrow g^*f^*T_Z.
$$
Because $f(\PP^1) \cap Y^\text{sm}$ is nonempty, the cokernel of $dg$
is torsion.  Because $f(\PP^1)$ is contained in a fiber of $g$,
$g^*f^* T_Z \cong \OO_{\PP^1_k}^r$ for some nonnegative integer $r$.
Because $f^*T_Y$ is generated by global sections, also the image of
$dg$ is generated by global sections.  But the only coherent subsheaf
of $\OO_{\PP^1_k}^r$ whose cokernel is torsion and that is generated
by global sections is \emph{all} of $\OO_{\PP^1_k}^r$.  Therefore $dg$
is surjective, i.e., $f(\PP^1) \subset Y^\text{sm}$.
\end{proof}

\medskip\noindent
By ~\cite[Thm. 1.1]{OS}, the Hilbert functor of $\overline{X}\rightarrow S$ is
represented by an algebraic space that is separated and locally
finitely presented, $\text{Hilb}_{\overline{X}/S}$.  And $W_i^\text{flat}
\subset X^\text{flat} \times_S \overline{X}$ is a closed substack that
is proper, flat and finitely presented
over $X^\text{flat}$.  Therefore there is a $1$-morphism of $S$-stacks,
$$
\phi^\text{flat}:X^\text{flat} \rightarrow \text{Hilb}_{\overline{X}/S},
$$
such that $W_i^\text{flat}$ is the pullback by $\phi^\text{flat}$ of
the universal closed substack.  
Denote by $X^*\subset X$ the maximal open substack over which
$\phi^\text{flat}$ extends to a morphism.  Denote by $Q^*\rightarrow
\text{Hilb}_{X/S}$ the Stein factorization of $X^*\rightarrow
\text{Hilb}_{X/S}$, i.e., the integral closure of the image in the
function field of the coarse moduli space $|X^*|$.   
Denote by $\phi:X^* \rightarrow Q^*$ the induced
morphism.  

\begin{prop} \label{prop-8}
\marpar{prop-8}
The morphism $\phi:X^* \rightarrow Q^*$ is a free rational quotient.
\end{prop}

\begin{proof}
By construction, $Q^*$ is normal, $Q^* \rightarrow S$ is
a finitely presented morphism whose geometric generic fiber is
integral, and $\phi$ is a dominant $1$-morphism whose geometric
generic fiber is integral.  It remains to prove Definition~\ref{defn-1} (ii).

\medskip\noindent
By ~\cite{OS}, $\text{Hilb}_{\overline{X}/S}$ satisfies the valuative criterion
of properness (but it is not necessarily proper since it is not
necessarily quasi-compact).  And $X$ is normal.  Therefore every
irreducible component of $X-X^*$ has codimension $\geq 2$ in $X$; more
precisely, the geometric generic fiber over $S$ has codimension $\geq
2$ in the geometric generic fiber of $X$ over $S$.  
For reasons similar to ~\cite[Prop.II.3.7]{K}, $X^*$ contains
$f(\PP^1)$ 
for every $H_j$ and general
$[f] \in H_j$.

\medskip\noindent
There is a dense open subspace $Q^0 \subset Q^*$ over which $\phi$
is flat and the geometric fibers are integral.  Replace $X^0$ by
$X^0\cap \phi^{-1}(Q^0)$.  
By Corollary~\ref{cor-5} (v), the subscheme $X^0\times Q^* X^0$ equals
$W_i \cap 
(X^0\times_S X^0)$.  In particular, for every $x\in X^0$ the fiber of
$\phi$ containing $x$ is $(W_i)_x$.  Therefore for a general fiber
of $\phi$, for a general pair of points in the fiber, there is a
free rational curve in $H_i$ whose image is contained in the fiber and
contains the two points.  By Lemma \ref{lem-7}, the image is contained
in the smooth locus of $\phi$, i.e. this is a free rational curve in
the fiber.  This proves Definition~\ref{defn-1} (ii).
\end{proof}

\begin{prop} \label{prop-9}
\marpar{prop-9}
The morphism $\phi:X^* \rightarrow Q^*$ is a maximal free rational
quotient.
\end{prop}

\begin{proof}
Let $\phi_1:X_1^* \rightarrow Q_1^*$ be a free rational quotient.  If
this is a trivial free rational quotient, the morphism $\psi$ is
trivial.  Therefore assume it is a nontrivial free rational quotient.

\medskip\noindent
There exists a dense open $U\subset Q_1^*$ such that,
$$
\chi: U' \rightarrow U,
$$
is faithfully flat and quasi-compact and the geometric fibers are integral, 
where $U' = X^0 \cap \phi_1^*(U)$ and where 
$\chi$ is the restriction of $\phi_1$.
By faithfully flat descent, to construct $\psi:U \rightarrow Q^*$, it
is equivalent to construct a morphism $\psi':U' \rightarrow Q^*$
satisfying a cocycle condition: indeed, 
the morphism $\psi$ is equivalent to the
graph of $\psi$, which is equivalent to a certain kind of
quasi-coherent sheaf on $U\times Q^*$, so faithfully flat descent
for quasi-coherent sheaves
applies to $(\chi,1):U'\times Q^* \rightarrow U\times Q^*$.  

\medskip\noindent
Define
$\psi'$ to be the restriction of $\phi$.  For each geometric point $x$
in $U'$, define $Y_x$ to be the fiber of $\chi$ containing $x$.  
The cocycle condition for $\psi'$ is that the fiber product $U'' = U'
\times_{\chi,U,\chi} U'$ is contained in $U' \times_{\psi',Q^*,\psi'}
U'$.  Now $U'' \rightarrow U'$ is a flat morphism whose geometric
fibers are integral.  Thus $U''$ is integral.  So it suffices to prove
it is set-theoretically contained in $W_i$, i.e., for a general geometric
point $x$ of $U'$, $U_x'' \subset (W_i)_x$.  

\medskip\noindent
By hypothesis, there is a
dense subset of $U_x''$ consisting of points $y$ contained in a free
morphism $f:\PP^1 \rightarrow X$ such that $f(0)=x$ and
$f(\infty)=y$.  By Corollary~\ref{cor-5} (ii), $f(\PP^1) \subset
(W_i)_x$, in particular $y\in (W_i)_x$.  Therefore $U_x'' \subset
(W_i)_x$.  
\end{proof}

\bibliography{my}
\bibliographystyle{abbrv}

\end{document}